\documentclass[11pt]{amsart}
\usepackage{amssymb,latexsym,amsfonts,verbatim,amscd,ifthen}
\usepackage[usenames,dvipsnames]{xcolor}
\usepackage[usenames,dvipsnames]{pstricks}
\usepackage{graphics,hyperref}
\usepackage{enumitem}
\usepackage{pstricks, pstricks-add}
\usepackage{epsfig}
\usepackage{pst-grad} 
\usepackage{pst-plot} 
\usepackage[normalem]{ulem}
\usepackage{algorithm}
\usepackage[noend]{algorithmic}

\definecolor{verylight}{gray}{0.50}
\definecolor{light}{gray}{0.9}
\definecolor{medium}{gray}{0.85}
\definecolor{orange}{cmy}{0,0.5,1}
\definecolor{orange}{rgb}{1,0.5,0}
\definecolor{orange}{RGB}{255,127,0}
\definecolor{orange}{cmyk}{0,0.5,1,0}
\newtheorem{thm1}{Theorem}[section]
\newtheorem{lem1}[thm1]{Lemma}

\newtheorem{rem1}[thm1]{Remark}

\newtheorem{def1}[thm1]{Definition}

\newtheorem{cor1}[thm1]{Corollary}

\newtheorem{prop1}[thm1]{Proposition}
\newtheorem{ex1}[thm1]{Example}

\DeclareMathOperator{\supp}{supp} 
\DeclareMathOperator{\chara}{char} 
\DeclareMathOperator{\Supp}{Supp}

  \DeclareMathOperator{\ini}{in}

\DeclareMathOperator{\Gr}{Gr}

\newcommand{\A}{\mathbb{A}}
\newcommand{\T}{\mathbb{T}}
\newcommand{\N}{\mathbb{N}}
\newcommand{\Z}{\mathbb{Z}}

\newcommand{\K}{\mathbb{K}} 
\newcommand{\C}{\mathbb{C}}
                                                                                                                                       
\begin{document}
\title{ Binomial fibers and indispensable binomials } 
\author[H. Charalambous]{Hara Charalambous}
\address{ Department of Mathematics, Aristotle University of Thessaloniki, Thessaloniki \\54124, Greece} \email{hara@math.auth.gr}
\author[A. Thoma]{ Apostolos Thoma}
\address{ Department of Mathematics, University of Ioannina, Ioannina 45110, Greece} \email{athoma@uoi.gr}
\author[M. Vladoiu]{Marius Vladoiu}
\address{ Faculty of Mathematics and Computer Science, University of Bucharest,  Str. Academiei 14, Bucharest,
RO-010014, Romania, and}
\address{Simion Stoilow Institute of Mathematics of Romanian Academy, Research group of the project PN-II-RU-TE-2012-3-0161, P.O.Box 1--764, Bucharest
014700, Romania}
\email{vladoiu@gta.math.unibuc.ro}

\thanks{The third author was partially supported by project PN-II-RU-TE-2012-3-0161, granted by the Romanian National Authority for Scientific Research, CNCS -
UEFISCDI. The paper was written when the third author visited University of Ioannina, and he gratefully acknowledges for the warm environment that fostered the collaboration of the authors.}
\keywords{Binomial ideals, Markov basis, lattices.}
\subjclass[2000]{13P10, 14M25}

\begin{abstract}
Let $I$ be an arbitrary ideal generated by binomials. We show that certain equivalence classes of fibers are associated to any minimal binomial generating set of $I$. We provide a simple and efficient algorithm to compute the indispensable binomials of a binomial ideal from a given generating set of binomials and an algorithm to detect whether a binomial ideal is generated by indispensable binomials.  
\end{abstract}
\maketitle

\section*{Introduction}
Let $R=\K[x_1,\ldots, x_n]$ where $\K$ is a field. A binomial is a polynomial of the form $x^{\bf u}- \lambda x^{\bf v}$ where ${\bf u},{\bf v}\in \N^n$  and $\lambda\in\K\setminus\{0\}$, and a binomial ideal is an ideal generated by binomials. We say that the ideal $I$ of $R$ is a {\it pure binomial ideal} if $I$ is generated by {\it pure difference binomials}, i.e. binomials of the form $x^{\bf u}-x^{\bf v}$ with ${\bf u},{\bf v}\in\N^n$. Binomial ideals were first studied systematically in \cite{ES} and this class of ideals also includes lattice ideals. Recall that if $L\subset \Z^n$ is a lattice, then the corresponding lattice ideal  is defined as $I_L=(x^{\bf u}-x^{\bf v}:\ {\bf u}-{\bf v}\in L)$ and the lattice is saturated exactly when the lattice ideal is toric, i.e. prime. The study of binomial ideals is a rich subject: the classical reference is  \cite{St} and we also refer to  \cite{EM} for recent developments. It has  applications in various areas in mathematics, such as algebraic  statistics,  integer programming, graph theory, computational biology,  code theory, see \cite{DS, DuZi, HT, OH, OH1, SS}, etc.

 A particular problem that arises  is  the {\it efficient} generation of binomial  ideals by a set of binomials. Up to now, it has mainly been addressed  for toric and lattice ideals, see \cite{BSR,  CKT, CTV, HM, HSt, St} among others. In Section 1, we consider this problem in the case of binomial ideals. For this we study the {\it fibers} of binomial ideals: in \cite[Proposition 2.4]{DMM} an equivalence relation  on $\N^n$ was introduced for any binomial ideal $I$ of $R$, namely ${\bf u}\sim_I {\bf v}$ if $x^{\bf u}- \lambda x^{\bf v}\in I$ for some $\lambda\neq 0$. For each such equivalence class, we get a {\it fiber} on the set of monomials: the $I$-fiber of $x^{\bf u}$ is   the set $ \{x^{\bf v}:\ \ {\bf u}\sim_I {\bf v}\}$. When $I:=I_L$ is the lattice ideal of $L$, the equivalence class of ${\bf u}$ consists precisely of all  ${\bf v}$ such that ${\bf u}-{\bf v} \in L$ and the $I$-fibers are finite exactly when $L\cap \N^n=\{ \bf 0\}$. In this case, for each $I$-fiber one can use a graph construction, see \cite{DS, CKT}, that determines the $I$-fibers that appear as invariants associated to  any minimal generating set of $I$. We also note that  in \cite{CTV} the fibers of $I_L$ were studied even when $L\cap \N^n\neq\{ \bf 0\}$.  In all cases  finite or not, it is clear that divisibility of monomials  does not induce necessarily a meaningful partial order on the set of $I$-fibers. In this paper for any binomial ideal $I$ we define an equivalence relation on the set of $I$-fibers and then order the equivalence classes of $I$-fibers, see Definition \ref{equiv}.  We note that  this was  first done in \cite{CTV} for the case $I=I_L$. This partial order allows us to prove that a certain set of equivalence classes of $I$-fibers is an invariant, associated to any generating set of $I$, see Theorem \ref{markov_fibers}. We note that lattice ideals have all fibers either finite or infinite and the equivalence classes of fibers for lattice ideals have the same cardinality, see \cite[Propositions 2.3 and 3.5]{CTV}. However this might not be the case for a binomial ideal and this constitutes an added degree of difficulty, see Examples~\ref{draw_fibers} c) and \ref{example_multi} c). Moreover we show that binomial ideals which contain monomials have a unique maximal fiber consisting of all monomials of $I$, see Theorem~\ref{unique_monomial_fiber}.  

A related question that attracted a lot of interest in the recent years is whether there is a unique  minimal binomial generating set for a binomial ideal. One of the first papers to deal with this question for lattice ideals from a purely theoretical point of view was \cite{PS}.  As it turns out, the positive answer has applications to Algebraic Statistics: \cite{AT, ATY, OH}. Thus in \cite{OH1} and \cite{ATY} the  notions of {\it indispensable} monomials and binomials were defined. Let $I$ be a binomial ideal. A binomial is called {\it indispensable} if (up to a nonzero constant) it belongs to every minimal generating set of $I$ consisting of binomials. This implies of course that (up to a nonzero constant) it belongs to every binomial generating set of $I$. A monomial is called {\it indispensable} if it is a monomial term of at least one binomial in every system of binomial generators of $I$. How does one compute these elements?

When $I:=I_L$ is a lattice ideal and $L\cap \N^n=\{ \bf 0\}$, there are several works in the literature that deal with this problem. In particular, in \cite{OH1} it was shown that to compute the indispensable binomials of $I_L$, one computes all lexicographic reduced Gr{\" o}bner bases and then their intersection: there are $n!$ such bases; a corresponding result for   indispensable monomials   was shown in \cite{ATY}. In \cite{OV},  it was shown that to compute the indispensable binomials of $I_L$, it is enough to compute certain degree-reverse lexicographic reduced Gr{\" o}bner bases of $I_L$ ($n$ of them), and then compute their intersection. In \cite[Proposition 3.1]{CKT}, it was shown that to find the indispensable monomials of $I_L$, it is enough to consider any one of the binomial generating sets of $I_L$. Moreover in \cite[Theorem 2.12]{CKT} it was shown that  in order to find the indispensable binomials of $I_L$, it is  enough to consider any  minimal binomial generating set of $I_L$, assign $\Z^n/L$-degrees to the  binomials of this set and to compute their minimal $\Z^n/L$-degrees. More recently in \cite[Theorem 1.1, Corollary 1.3]{KO}, it was shown  that if $I$ is a pure binomial ideal then there is a $d\in \N$ such that any  $I$ is $A$-graded for some $A\subset \Z^d$:  when  $\N A\cap (-\N A)=\{{\bf 0}\}$ and all fibers are finite a sufficient condition was given in \cite{KO} for the indispensable binomials and a characterization for the indispensable monomials, involving the $I$-fibers of a minimal generating set of $I$.  

In this paper, we significantly improve all previously known results  regarding indispensable binomials. Moreover  our  results apply  to the general case of all  binomial ideals. In Section 2 we show that as in \cite{KO}, the indispensable monomials are the elements of the minimal generating set of the monomial ideal of $I$, see Remark~\ref{known_indisp}. Then we go on and are able to express this condition into three necessary and equivalent conditions involving  a graph whose vertices are the (possibly infinitely many) elements of the fiber, see Theorem~\ref{indisp_mon_fiber}.  This result is then applied to provide sufficient and necessary conditions for a binomial in $I$ to be indispensable, see Theorem~\ref{indisp_bin_fiber}.

In Section 3, we prove that an arbitrary system of binomial generators of $I$ gives all necessary information to decide whether a given binomial is indispensable,  see Theorem \ref{algorithm_indisp_thm}. As an immediate application of Theorem~\ref{algorithm_indisp_thm} we obtain an algorithm which computes the indispensable elements of a binomial ideal $I$, given any system of binomial generators of $I$, see Algorithm~\ref{indisp_binomials}. This algorithm bypasses the computation of a reduced Gr\" obner basis unlike the previous methods. As a result we show that Algorithm~\ref{indisp_binomials} is a polynomial time algorithm, see Remark~\ref{polynomial_time_complexity}, in contrast to the other methods, see \cite{Ma} for details regarding the complexity of reduced Gr\" obner basis computation for binomial ideals. We also show, that to decide whether a minimal system of binomial generators of a binomial ideal is in fact a system of indispensable binomials it is enough to compute the cardinality of the minimal generating set of an associated monomial ideal, see Corollary~\ref{gen_indispensable_cor} and the resulting Algorithm \ref{gen_indispensable}.  
 
In Section 4, we generalize the notion of primitive elements to pure binomial ideals. The set of all primitive elements is the Graver basis of $I$. This set is extremely important in theory and all computations involving lattice ideals, see \cite{St}. We prove that the Graver basis of any pure binomial ideal is finite, see Proposition~\ref{finit_grav}. We show that the Graver basis includes as a subset the universal Groebner basis of I, see Theorem \ref{ugb_gr}. Finally, we show that a Lawrence lifting construction gives a pure binomial ideal generated by indispensable binomials, see Theorem \ref{arbitrary_indispensable}.  

We thank the two anonymous referees for their careful reading and suggestions which greatly improved our paper, in particular for the questions on the complexity of the algorithms and on the generalization to binomial ideals, which we could answer affirmatively.

\section{Fibers of a Binomial Ideal}

Let $R=\K[x_1,\ldots, x_n]$ where $\K$ is a field. We denote by $\T^n$ the set of monomials of $R$ including $1=x^0$, where $x^{\bf u}=x_1^{u_1}\cdots x_n^{u_n}$. If $J$ is a monomial ideal of $R$ we denote by $G(J)$ the unique set of minimal monomial generators of $J$. For $B=x^{\bf u}-\lambda x^{\bf v}$ with $\lambda\neq 0$ we let $\supp(B):=\{x^{\bf u},x^{\bf v}\}$. 

\begin{def1}
{\em Let $I$ be a binomial ideal of $R$. We say that $F\subset \T^n$ is an $I$-{\it  fiber} if there exists a $x^{\bf u}\in \T^n$ such that $F= \{ x^{\bf v}:\ {\bf v}\sim_I {\bf u} \}$. If $x^{\bf u}\in \T^n$, and  $F$ is an $I$-fiber containing $x^{\bf u}$ we write $F_{\bf u}$ or $F_{x^{\bf u}}$ for $F$. If  $B\in I$ and $B=x^{\bf u}-\lambda x^{\bf v}$ with $\lambda\neq 0$ we write $F_B$ for $F_{\bf u}$.}
\end{def1}

It is trivial that $|F_{\bf u}|=1$, that is $F_{\bf u}$ is a singleton, if and only if there is no binomial $0\neq B\in I$ such that $x^{\bf u}\in \supp(B)$.
If $J\subset I$ gives the containment between two binomial ideals and $F$ is a $J$-fiber, then clearly $F$ is contained in an $I$-fiber.

\begin{ex1}\label{draw_fibers}$ $
\newline
\noindent
{\em a) Let $a\in\N$,  $r\in\Z_{\geq 1}$, $L=r\Z$, $\lambda\in\K\setminus\{0\}$ and $I_1=(x^a-\lambda x^{a+r})$. It is immediate that $x^{a+j}-\lambda^{nr}x^{a+nr+j}\in I_1$ for all $j,n\in \N$. The $I_1$-fibers are either singletons or infinite. There are exactly $a$ singletons and $r$ distinct infinite $I_1$-fibers  of the form  $F_{a+j}=\{x^{a+j+nr}:\ n\in\N\}$, for  $0\leq j\leq r-1$. Moreover 
since 
\[
x^a-\lambda x^{a+r}=(1-\lambda x^r)(x^a-\lambda^{3r}x^{a+3r})+\lambda^{3r}x^{3r-1}(x^{a+1}-\lambda x^{a+r+1})\ ,
\]
$I_1=(x^a-\lambda^{3r}x^{a+3r}, x^a-\lambda x^{a+ r+1})$. Thus $I_1$ has no indispensable binomials. It is  clear that the the only indispensable monomial of $I_1$ is $x^a$. 

\noindent b) Let $I_2=(y-x^2y,y^3-xy^3,y^4-9y^6,y^7-3y^8)$.  The $I_2$-fibers are as follows:
\begin{itemize}
\item{} $F_{x^i}=\{x^i\}$ for all $i\in\N$
\item{} $F_y=\{yx^{2n}:\  n\in \N \}$
\item{} $F_{yx}=\{ yx^{2n+1}:\ n\in\N\}$
\item{} $F_{y^2}=\{y^2x^{2n}:\ n\in\N\}$
\item{} $F_{y^2x}=\{ y^2x^{2n+1}:\ n\in\N\}$
\item{} $F_{y^3}=\{y^3x^{n}:\ n\in\N\}$
\item{} $F_{y^4}=\{y^{4+m}x^{n}:\ m,n\in\N\}$,
\end{itemize}
and they are independent of the characteristic of $\K$. Note that when $\chara(\K)=3$ then $I_2=(y-x^2y,y^3-xy^3,y^4)$ and is the binomial ideal $(y-x^2y,y^3-xy^3,y^4-y^5,y^4+y^5)$. The $I_2$-fibers are depicted in the left of Figure~\ref{fibers}.

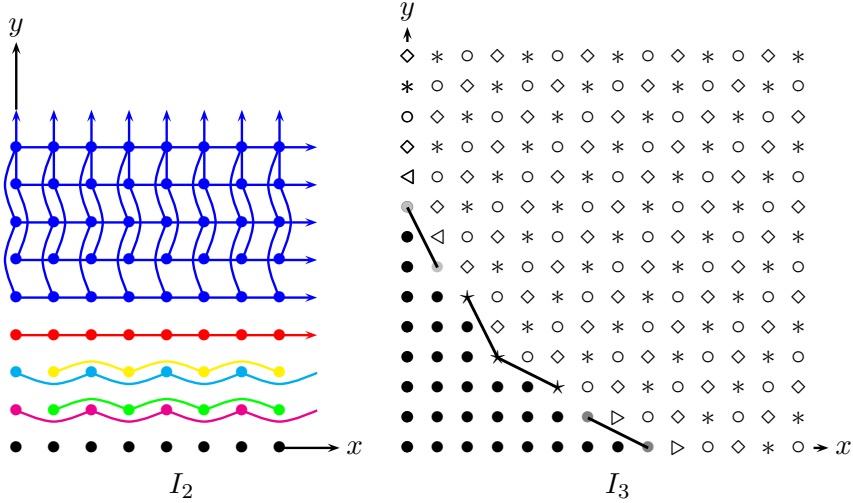
\begin{figure}[hbt]
\begin{center}
\psset{unit=1cm} 
\begin{pspicture}(0,-0.5)(11,5.8)
\rput(0,0){$\bullet$} \rput(0.5,0){$\bullet$} \rput(1,0){$\bullet$} \rput(1.5,0){$\bullet$} \rput(2,0){$\bullet$} \rput(2.5,0){$\bullet$} \rput(3,0){$\bullet$} \rput(3.5,0){$\bullet$} 
\rput(0,0.5){\magenta $\bullet$} \rput(0.5,0.5){\green $\bullet$} \rput(1,0.5){\magenta $\bullet$} \rput(1.5,0.5){\green $\bullet$} \rput(2,0.5){\magenta $\bullet$} \rput(2.5,0.5){\green $\bullet$} \rput(3,0.5){\magenta $\bullet$} \rput(3.5,0.5){\green $\bullet$}
\rput(0,1){\cyan $\bullet$} \rput(0.5,1){\yellow $\bullet$} \rput(1,1){\cyan $\bullet$} \rput(1.5,1){\yellow $\bullet$} \rput(2,1){\cyan $\bullet$} \rput(2.5,1){\yellow $\bullet$} \rput(3,1){\cyan $\bullet$} \rput(3.5,1){\yellow $\bullet$}
\rput(0,1.5){\red $\bullet$} \rput(0.5,1.5){\red $\bullet$} \rput(1,1.5){\red $\bullet$} \rput(1.5,1.5){\red $\bullet$} \rput(2,1.5){\red $\bullet$} \rput(2.5,1.5){\red $\bullet$} \rput(3,1.5){\red $\bullet$} \rput(3.5,1.5){\red $\bullet$}
\rput(0,2){\blue $\bullet$} \rput(0.5,2){\blue $\bullet$} \rput(1,2){\blue $\bullet$} \rput(1.5,2){\blue $\bullet$} \rput(2,2){\blue $\bullet$} \rput(2.5,2){\blue $\bullet$} \rput(3,2){\blue $\bullet$} \rput(3.5,2){\blue $\bullet$}
\rput(0,2.5){\blue $\bullet$} \rput(0.5,2.5){\blue $\bullet$} \rput(1,2.5){\blue $\bullet$} \rput(1.5,2.5){\blue $\bullet$} \rput(2,2.5){\blue $\bullet$} \rput(2.5,2.5){\blue $\bullet$} \rput(3,2.5){\blue $\bullet$} \rput(3.5,2.5){\blue $\bullet$}
\rput(0,3){\blue $\bullet$} \rput(0.5,3){\blue $\bullet$} \rput(1,3){\blue $\bullet$} \rput(1.5,3){\blue $\bullet$} \rput(2,3){\blue $\bullet$} \rput(2.5,3){\blue $\bullet$} \rput(3,3){\blue $\bullet$} \rput(3.5,3){\blue $\bullet$}
\rput(0,3.5){\blue $\bullet$} \rput(0.5,3.5){\blue $\bullet$} \rput(1,3.5){\blue $\bullet$} \rput(1.5,3.5){\blue $\bullet$} \rput(2,3.5){\blue $\bullet$} \rput(2.5,3.5){\blue $\bullet$} \rput(3,3.5){\blue $\bullet$} \rput(3.5,3.5){\blue $\bullet$}
\rput(0,4){\blue $\bullet$} \rput(0.5,4){\blue $\bullet$} \rput(1,4){\blue $\bullet$} \rput(1.5,4){\blue $\bullet$} \rput(2,4){\blue $\bullet$} \rput(2.5,4){\blue $\bullet$} \rput(3,4){\blue $\bullet$} \rput(3.5,4){\blue $\bullet$}
\pscurve[linewidth=0.9pt, linecolor=magenta](0,0.5)(0.5,0.35)(1,0.5)
\pscurve[linewidth=0.9pt, linecolor=magenta](1,0.5)(1.5,0.35)(2,0.5)
\pscurve[linewidth=0.9pt, linecolor=magenta](2,0.5)(2.5,0.35)(3,0.5)
\pscurve[linewidth=0.9pt, linecolor=magenta](3,0.5)(3.5,0.35)(4,0.5)

\pscurve[linewidth=0.9pt, linecolor=green](0.5,0.5)(1,0.65)(1.5,0.5)
\pscurve[linewidth=0.9pt, linecolor=green](1.5,0.5)(2,0.65)(2.5,0.5)
\pscurve[linewidth=0.9pt, linecolor=green](2.5,0.5)(3,0.65)(3.5,0.5)

\pscurve[linewidth=0.9pt, linecolor=cyan](0,1)(0.5,0.85)(1,1)
\pscurve[linewidth=0.9pt, linecolor=cyan](1,1)(1.5,0.85)(2,1)
\pscurve[linewidth=0.9pt, linecolor=cyan](2,1)(2.5,0.85)(3,1)
\pscurve[linewidth=0.9pt, linecolor=cyan](3,1)(3.5,0.85)(4,1)

\pscurve[linewidth=0.9pt, linecolor=yellow](0.5,1)(1,1.15)(1.5,1)
\pscurve[linewidth=0.9pt, linecolor=yellow](1.5,1)(2,1.15)(2.5,1)
\pscurve[linewidth=0.9pt, linecolor=yellow](2.5,1)(3,1.15)(3.5,1)

\psline[linewidth=0.9pt, linecolor=red](0,1.5)(3.5,1.5)

\psline[linewidth=0.9pt, linecolor=blue](0,2)(3.5,2)
\psline[linewidth=0.9pt, linecolor=blue](0,3)(3.5,3)
\pscurve[linewidth=0.9pt, linecolor=blue](0,2)(-0.15,2.5)(0,3)
\pscurve[linewidth=0.9pt, linecolor=blue](0.5,2)(0.35,2.5)(0.5,3)
\pscurve[linewidth=0.9pt, linecolor=blue](1,2)(0.85,2.5)(1,3)
\pscurve[linewidth=0.9pt, linecolor=blue](1.5,2)(1.35,2.5)(1.5,3)
\pscurve[linewidth=0.9pt, linecolor=blue](2,2)(1.85,2.5)(2,3)
\pscurve[linewidth=0.9pt, linecolor=blue](2.5,2)(2.35,2.5)(2.5,3)
\pscurve[linewidth=0.9pt, linecolor=blue](3,2)(2.85,2.5)(3,3)
\pscurve[linewidth=0.9pt, linecolor=blue](3.5,2)(3.35,2.5)(3.5,3)
\pscurve[linewidth=0.9pt, linecolor=blue](0,3)(-0.15,3.5)(0,4)
\pscurve[linewidth=0.9pt, linecolor=blue](0.5,3)(0.35,3.5)(0.5,4)
\pscurve[linewidth=0.9pt, linecolor=blue](1,3)(0.85,3.5)(1,4)
\pscurve[linewidth=0.9pt, linecolor=blue](1.5,3)(1.35,3.5)(1.5,4)
\pscurve[linewidth=0.9pt, linecolor=blue](2,3)(1.85,3.5)(2,4)
\pscurve[linewidth=0.9pt, linecolor=blue](2.5,3)(2.35,3.5)(2.5,4)
\pscurve[linewidth=0.9pt, linecolor=blue](3,3)(2.85,3.5)(3,4)
\pscurve[linewidth=0.9pt, linecolor=blue](3.5,3)(3.35,3.5)(3.5,4)

\psline[linewidth=0.9pt, linecolor=blue](0,2.5)(3.5,2.5)
\psline[linewidth=0.9pt, linecolor=blue](0,3.5)(3.5,3.5)
\psline[linewidth=0.9pt, linecolor=blue](0,4)(3.5,4)
\psline[linewidth=0.9pt, linecolor=blue](0,3.5)(0,4)
\psline[linewidth=0.9pt, linecolor=blue](0.5,3.5)(0.5,4)
\psline[linewidth=0.9pt, linecolor=blue](1,3.5)(1,4)
\psline[linewidth=0.9pt, linecolor=blue](1.5,3.5)(1.5,4)
\psline[linewidth=0.9pt, linecolor=blue](2,3.5)(2,4)
\psline[linewidth=0.9pt, linecolor=blue](2.5,3.5)(2.5,4)
\psline[linewidth=0.9pt, linecolor=blue](3,3.5)(3,4)
\psline[linewidth=0.9pt, linecolor=blue](3.5,3.5)(3.5,4)
\pscurve[linewidth=0.9pt, linecolor=blue](0,2.5)(0.15,3)(0,3.5)
\pscurve[linewidth=0.9pt, linecolor=blue](0.5,2.5)(0.65,3)(0.5,3.5)
\pscurve[linewidth=0.9pt, linecolor=blue](1,2.5)(1.15,3)(1,3.5)
\pscurve[linewidth=0.9pt, linecolor=blue](1.5,2.5)(1.65,3)(1.5,3.5)
\pscurve[linewidth=0.9pt, linecolor=blue](2,2.5)(2.15,3)(2,3.5)
\pscurve[linewidth=0.9pt, linecolor=blue](2.5,2.5)(2.65,3)(2.5,3.5)
\pscurve[linewidth=0.9pt, linecolor=blue](3,2.5)(3.15,3)(3,3.5)
\pscurve[linewidth=0.9pt, linecolor=blue](3.5,2.5)(3.65,3)(3.5,3.5)


\psline[linewidth=0.9pt, linecolor=red]{->}(3.5,1.5)(4,1.5)
\psline[linewidth=0.9pt, linecolor=blue]{->}(3.5,2)(4,2)
\psline[linewidth=0.9pt, linecolor=blue]{->}(3.5,3)(4,3)
\psline[linewidth=0.9pt, linecolor=blue]{->}(3.5,2.5)(4,2.5)
\psline[linewidth=0.9pt, linecolor=blue]{->}(3.5,3.5)(4,3.5)
\psline[linewidth=0.9pt, linecolor=blue]{->}(3.5,4)(4,4)
\psline[linewidth=0.9pt, linecolor=blue]{->}(0,4)(0,4.5)
\psline[linewidth=0.9pt, linecolor=blue]{->}(0.5,4)(0.5,4.5)
\psline[linewidth=0.9pt, linecolor=blue]{->}(1,4)(1,4.5)
\psline[linewidth=0.9pt, linecolor=blue]{->}(1.5,4)(1.5,4.5) 
\psline[linewidth=0.9pt, linecolor=blue]{->}(2,4)(2,4.5)
\psline[linewidth=0.9pt, linecolor=blue]{->}(2.5,4)(2.5,4.5)
\psline[linewidth=0.9pt, linecolor=blue]{->}(3,4)(3,4.5)
\psline[linewidth=0.9pt, linecolor=blue]{->}(3.5,4)(3.5,4.5)

\psline[linewidth=0.9pt, linecolor=black]{->}(3.5,0)(4.3,0)
\psline[linewidth=0.9pt, linecolor=black]{->}(0,4.5)(0,5.4)
\rput(4.5,0){\bf $x$}
\rput(0,5.6){\bf $y$}

\rput(5.2,0){$\bullet$} \rput(5.6,0){$\bullet$} \rput(6,0){$\bullet$} \rput(6.4,0){$\bullet$} \rput(6.8,0){$\bullet$} \rput(7.2,0){$\bullet$} \rput(7.6,0){$\bullet$} \rput(8,0){$\bullet$} \rput(8.4,0){\gray $\bullet$} \rput(8.8,0){$\triangleright$} \rput(9.2,0){$\circ$} \rput(9.6,0){$\diamond$} \rput(10,0){$\ast$} \rput(10.4,0){$\circ$}
\psline[linewidth=0.9pt, linecolor=black]{->}(10.6,0)(10.8,0)
\rput(11,0){\bf $x$}

\rput(5.2,0){$\bullet$} \rput(5.2,0.4){$\bullet$} \rput(5.2,0.8){$\bullet$} \rput(5.2,1.2){$\bullet$} \rput(5.2,1.6){$\bullet$} \rput(5.2,2){$\bullet$} \rput(5.2,2.4){$\bullet$} \rput(5.2,2.8){$\bullet$} \rput(5.2,3.2){$\bullet$} \rput(5.2,3.6){$\triangleleft$} \rput(5.2,4){$\diamond$} \rput(5.2,4.4){$\circ$} \rput(5.2,4.8){$\ast$} \rput(5.2,5.2){$\diamond$}
\psline[linewidth=0.9pt, linecolor=black]{->}(5.2,5.4)(5.2,5.6)
\rput(5.2,5.8){\bf $y$}

\rput(5.2,0.4){$\bullet$} \rput(5.6,0.4){$\bullet$} \rput(6,0.4){$\bullet$} \rput(6.4,0.4){$\bullet$} \rput(6.8,0.4){$\bullet$} \rput(7.2,0.4){$\bullet$} \rput(7.6,0.4){\gray $\bullet$} \rput(8,0.4){$\triangleright$} \rput(8.4,0.4){$\circ$} \rput(8.8,0.4){$\diamond$} \rput(9.2,0.4){$\ast$} \rput(9.6,0.4){$\circ$} \rput(10,0.4){$\diamond$} \rput(10.4,0.4){$\ast$}

\rput(5.2,0.8){$\bullet$} \rput(5.6,0.8){$\bullet$} \rput(6,0.8){$\bullet$} \rput(6.4,0.8){$\bullet$} \rput(6.8,0.8){$\bullet$} \rput(7.2,0.8){$\star$} \rput(7.6,0.8){$\circ$} \rput(8,0.8){$\diamond$} \rput(8.4,0.8){$\ast$} \rput(8.8,0.8){$\circ$} \rput(9.2,0.8){$\diamond$} \rput(9.6,0.8){$\ast$} \rput(10,0.8){$\circ$} \rput(10.4,0.8){$\diamond$}

\rput(5.2,1.2){$\bullet$} \rput(5.6,1.2){$\bullet$} \rput(6,1.2){$\bullet$} \rput(6.4,1.2){$\star$} \rput(6.8,1.2){$\circ$} \rput(7.2,1.2){$\diamond$} \rput(7.6,1.2){$\ast$} \rput(8,1.2){$\circ$} \rput(8.4,1.2){$\diamond$} \rput(8.8,1.2){$\ast$} \rput(9.2,1.2){$\circ$} \rput(9.6,1.2){$\diamond$} \rput(10,1.2){$\ast$} \rput(10.4,1.2){$\circ$}

\rput(5.2,1.6){$\bullet$} \rput(5.6,1.6){$\bullet$} \rput(6,1.6){$\bullet$} \rput(6.4,1.6){$\diamond$} \rput(6.8,1.6){$\ast$} \rput(7.2,1.6){$\circ$} \rput(7.6,1.6){$\diamond$} \rput(8,1.6){$\ast$} \rput(8.4,1.6){$\circ$} \rput(8.8,1.6){$\diamond$} \rput(9.2,1.6){$\ast$} \rput(9.6,1.6){$\circ$} \rput(10,1.6){$\diamond$} \rput(10.4,1.6){$\ast$}

\rput(5.2,2){$\bullet$} \rput(5.6,2){$\bullet$} \rput(6,2){$\star$} \rput(6.4,2){$\circ$} \rput(6.8,2){$\diamond$} \rput(7.2,2){$\ast$} \rput(7.6,2){$\circ$} \rput(8,2){$\diamond$} \rput(8.4,2){$\ast$} \rput(8.8,2){$\circ$} \rput(9.2,2){$\diamond$} \rput(9.6,2){$\ast$} \rput(10,2){$\circ$} \rput(10.4,2){$\diamond$}

\rput(5.2,2.4){$\bullet$} \rput(5.6,2.4){\lightgray $\bullet$} \rput(6,2.4){$\diamond$} \rput(6.4,2.4){$\ast$} \rput(6.8,2.4){$\circ$} \rput(7.2,2.4){$\diamond$} \rput(7.6,2.4){$\ast$} \rput(8,2.4){$\circ$} \rput(8.4,2.4){$\diamond$} \rput(8.8,2.4){$\ast$} \rput(9.2,2.4){$\circ$} \rput(9.6,2.4){$\diamond$} \rput(10,2.4){$\ast$} \rput(10.4,2.4){$\circ$}

\rput(5.2,2.8){$\bullet$} \rput(5.6,2.8){$\triangleleft$} \rput(6,2.8){$\circ$} \rput(6.4,2.8){$\diamond$} \rput(6.8,2.8){$\ast$} \rput(7.2,2.8){$\circ$} \rput(7.6,2.8){$\diamond$} \rput(8,2.8){$\ast$} \rput(8.4,2.8){$\circ$} \rput(8.8,2.8){$\diamond$} \rput(9.2,2.8){$\ast$} \rput(9.6,2.8){$\circ$} \rput(10,2.8){$\diamond$} \rput(10.4,2.8){$\ast$}

\rput(5.2,3.2){\lightgray $\bullet$} \rput(5.6,3.2){$\diamond$} \rput(6,3.2){$\ast$} \rput(6.4,3.2){$\circ$} \rput(6.8,3.2){$\diamond$} \rput(7.2,3.2){$\ast$} \rput(7.6,3.2){$\circ$} \rput(8,3.2){$\diamond$} \rput(8.4,3.2){$\ast$} \rput(8.8,3.2){$\circ$} \rput(9.2,3.2){$\diamond$} \rput(9.6,3.2){$\ast$} \rput(10,3.2){$\circ$} \rput(10.4,3.2){$\diamond$}

\rput(5.2,3.6){$\triangleleft$} \rput(5.6,3.6){$\circ$} \rput(6,3.6){$\diamond$} \rput(6.4,3.6){$\ast$} \rput(6.8,3.6){$\circ$} \rput(7.2,3.6){$\diamond$} \rput(7.6,3.6){$\ast$} \rput(8,3.6){$\circ$} \rput(8.4,3.6){$\diamond$} \rput(8.8,3.6){$\ast$} \rput(9.2,3.6){$\circ$} \rput(9.6,3.6){$\diamond$} \rput(10,3.6){$\ast$} \rput(10.4,3.6){$\circ$}

\rput(5.2,4){$\diamond$} \rput(5.6,4){$\ast$} \rput(6,4){$\circ$} \rput(6.4,4){$\diamond$} \rput(6.8,4){$\ast$} \rput(7.2,4){$\circ$} \rput(7.6,4){$\diamond$} \rput(8,4){$\ast$} \rput(8.4,4){$\circ$} \rput(8.8,4){$\diamond$} \rput(9.2,4){$\ast$} \rput(9.6,4){$\circ$} \rput(10,4){$\diamond$} \rput(10.4,4){$\ast$}

\rput(5.2,4.4){$\circ$} \rput(5.6,4.4){$\diamond$} \rput(6,4.4){$\ast$} \rput(6.4,4.4){$\circ$} \rput(6.8,4.4){$\diamond$} \rput(7.2,4.4){$\ast$} \rput(7.6,4.4){$\circ$} \rput(8,4.4){$\diamond$} \rput(8.4,4.4){$\ast$} \rput(8.8,4.4){$\circ$} \rput(9.2,4.4){$\diamond$} \rput(9.6,4.4){$\ast$} \rput(10,4.4){$\circ$} \rput(10.4,4.4){$\diamond$}

\rput(5.2,4.8){$\ast$} \rput(5.6,4.8){$\circ$} \rput(6,4.8){$\diamond$} \rput(6.4,4.8){$\ast$} \rput(6.8,4.8){$\circ$} \rput(7.2,4.8){$\diamond$} \rput(7.6,4.8){$\ast$} \rput(8,4.8){$\circ$} \rput(8.4,4.8){$\diamond$} \rput(8.8,4.8){$\ast$} \rput(9.2,4.8){$\circ$} \rput(9.6,4.8){$\diamond$} \rput(10,4.8){$\ast$} \rput(10.4,4.8){$\circ$}

\rput(5.2,5.2){$\diamond$} \rput(5.6,5.2){$\ast$} \rput(6,5.2){$\circ$} \rput(6.4,5.2){$\diamond$} \rput(6.8,5.2){$\ast$} \rput(7.2,5.2){$\circ$} \rput(7.6,5.2){$\diamond$} \rput(8,5.2){$\ast$} \rput(8.4,5.2){$\circ$} \rput(8.8,5.2){$\diamond$} \rput(9.2,5.2){$\ast$} \rput(9.6,5.2){$\circ$}\rput(10,5.2){$\diamond$} \rput(10.4,5.2){$\ast$}

\psline[linewidth=1.1pt](5.2,3.2)(5.6,2.4)
\psline[linewidth=1.1pt](6,2)(6.4,1.2)
\psline[linewidth=1.1pt](6.4,1.2)(7.2,0.8)
\psline[linewidth=1.1pt](7.6,0.4)(8.4,0)

\rput(2.2,-0.5){$I_2$} 
\rput(8,-0.5){$I_3$}

\end{pspicture}
\end{center}
\caption{Fibers of binomial ideals}\label{fibers}
\end{figure}

\noindent c) Consider the pure binomial ideal $I_3=(y^8-xy^6,x^2y^5-x^3y^3,x^3y^3-x^5y^2,x^6y-x^8)$. Its fibers are depicted in the right part of Figure~\ref{fibers}. There are 29  singleton $I$-fibers,  depicted by dots. The other fibers are:
 \begin{itemize}
 \item{} $F_{xy^6}=\{xy^6,y^8\}$ 
 \item{} $F_{x^3y^3}=\{x^3y^3,x^2y^5,x^5y^2\}$
 \item{} $F_{x^6y}=\{x^6y,x^8\}$
  \item{} $F_{xy^7}=\{xy^7,y^9\}$
  \item{} $F_{x^7y}=\{x^7y,x^9\}$
  \item{} $F_{x^4y^3}=\{ x^ay^b| \ (a,b)\in\N^2 \text{ and } (a,b)\in (4,3)+\N (-1,2)+\N (2,-1)\}$
  \item{} $F_{x^3y^4}=\{ x^ay^b| \ (a,b)\in\N^2 \text{ and } (a,b)\in (3,4)+\N (-1,2)+\N (2,-1)\}$
  \item{} $F_{x^4y^4}=\{ x^ay^b| \ (a,b)\in\N^2 \text{ and } (a,b)\in (4,4)+\N (-1,2)+\N (2,-1)\}$.
 \end{itemize}

\noindent d) Consider the ideal $I_4=(y^8-xy^6,x^2y^5-x^3y^3,x^3y^3-x^5y^2,x^6y-5x^8)\subset\K[x,y]$. If $\K$ is a field of characteristic different from 5, then all of the fibers are exactly the same as in the previous example (c) except of the three fibers $F_{x^4y^3}, F_{x^3y^4}, F_{x^4y^4}$ of $I_3$ which become one fiber of $I_4$. On the other hand, if $\chara(\K)=5$ then $I_4=(y^8-xy^6,x^2y^5-x^3y^3,x^3y^3-x^5y^2,x^6y)$ is the binomial ideal $(y^8-xy^6,x^2y^5-x^3y^3,x^3y^3-x^5y^2,x^6y-x^4y^3,x^6y+x^4y^3)$ and there are infinitely many singleton $I_4$-fibers: all of the singleton fibers of example (c) and all of the singleton fibers corresponding to monomials of the form $x^n$. The other $I_4$-fibers are:
\begin{itemize}
\item{} $F_{xy^6}=\{xy^6,y^8\}$ 
\item{} $F_{x^3y^3}=\{x^3y^3,x^2y^5,x^5y^2\}$
\item{} $F_{xy^7}=\{xy^7,y^9\}$
\item{} $F_{x^6y}$ which contains all the monomials belonging to the monomial ideal $M=(x^6y, x^4y^3, x^3y^4, x^2y^6, xy^8, y^{10})$.
\end{itemize}

Note that this example reveals that fibers are dependent on the characteristic of the field.
}
\end{ex1}

If $G\subset \mathbb{T}^n$ and ${\bf t}\in \N^n$ we let $x^{\bf t} G:=\{ x^{\bf t} x^{\bf u}:\  x^{\bf u}\in G\}$.

\begin{thm1} \label{partition_classif}
Let $\mathbf{F}$ be a partition of $\mathbb{T}^n$. There exists a binomial  ideal $I$ such that $\mathbf{F}$ is the set of $I$-fibers if and only if for any $ {\bf u}\in \mathbb{N}^n$ and any  $F\in\mathbf{F}$ there exists a $G\in \mathbf{F}$ such that $x^{\bf u}F\subset G$.
\end{thm1}

\begin{proof} Let  $I$ be a binomial ideal. Let $F:=F_{\bf t}$ be an $I$-fiber and let ${\bf u}\in \mathbb{N}^n$. It is clear  $x^{\bf u}F\subset F_{{\bf t}+{\bf u}}$. For the converse we let ${\bf F}=(F_i:\ i\in\Lambda)$ and $I$ be the ideal generated by the set $\{ x^{\bf u}-x^{\bf v}:\ {\bf u},{\bf v}\in F_i, i\in\Lambda\}$. Consider the set  ${\bf G}$ of $I$-fibers.  We will show that $\mathbf{G}=\mathbf{F}$. Indeed let $G$ be an $I$-fiber and $x^{\bf u}\in G$. Since $\mathbf{F}$ is a partition of $\mathbb{T}^n$, there is an $F\in \mathbf{F}$ such that $x^{\bf u}\in F$. From the definition of $I$, it is clear that $F\subset G$. For the converse inclusion suppose that $x^{\bf v}\in G$ and thus   $x^{\bf v}-x^{\bf u}\in I$. Hence $x^{\bf v}-x^{\bf u}=\sum_{k=1}^l c_k x^{{\bf w}_k}(x^{{\bf v}_k}-x^{{\bf u}_k})$ where $x^{{\bf u}_k},x^{{\bf v}_k}\in F_{i_k}$, $i_k\in\Lambda$, $c_k\in\K$, and $l\geq 1$. Therefore for some $k$, $x^{{\bf w}_k}x^{{\bf u}_k}=x^{\bf u}\in F$. By the hypothesis on the elements of $\mathbf{F}$ we obtain $x^{{\bf w}_k}x^{{\bf v}_k}\in F$. An easy induction on $l$ finishes the proof.
\end{proof}

In the proof of Theorem~\ref{partition_classif} we showed even more, namely that for a given set of fibers there exists a pure difference binomial ideal whose set of $I$-fibers is the given one.

A vector ${\bf u}\in \Z^n$  is called {\it pure} if ${\bf u}\in\N^n$ or $-{\bf u}\in\N^n$. If $B=x^{\bf u}-\lambda x^{\bf v}\in I$ with $\lambda\neq 0$ we let ${\bf v}(B)= {\bf u}-{\bf v}$. If $F$ is a fiber of a binomial ideal $I$, we let $L_F:=\langle {\bf v}(B):\  \ F_B=F\rangle\subset \Z^n$. We also consider $L_{pure,F}$, the sublattice of $L_F$ generated by the set
\[
\{{\bf w}\in\N^n:\ \exists \ x^{\bf u},x^{\bf v}\in F \text{ such that } {\bf w}={\bf u}-{\bf v}\},
\] 
and denote by $L_{pure, F}^+ $ the semigroup generated by the same set. Finally we let $M_F$ be the monomial ideal of $R$ generated by the elements of $F$.  In some cases it  might be that the monomials in $M_F$ are precisely the elements of $F$, for example  if $I=(1-x)$ in $\K[x]$, but usually this is far from being the case. If $I$ is a lattice ideal and $I=I_L$ then $L_F \subset L$ and by \cite[Proposition 2.3]{CTV} or \cite[Theorem 8.6]{MS} it follows that $F$ is an infinite fiber if and only if  $L_F$ contains a nonzero pure element. In the case of an arbitrary binomial ideal $I$, one can extend Proposition 2.6 of  \cite{CTV}.

\begin{prop1} Let $I$ be a binomial ideal and $F$ be an $I$-fiber.  If $\{x^{{\bf a}_1},\ldots,x^{{\bf a}_s}\}$ is the minimal monomial generating set of $M_F$ then
\[
F=\bigcup^s_{i=1}\{x^{{\bf a}_i}x^{\bf w}:\ {\bf w}\in L^{+}_{pure,F}\}\ .
\]
In particular, $F$ is infinite if and only if $L_{pure,F}\neq 0$.
\end{prop1}
\begin{proof}
It follows immediately from definitions the inclusion $\subseteq$. For the other inclusion let ${\bf w}\in L^+_{pure,F}$ and fix an $i$. Then there exist $x^{\bf u}\in F$ and $\lambda\in\K\setminus\{0\}$ such that $x^{\bf u}-\lambda x^{{\bf u}+{\bf w}}\in I$. Since $x^{{\bf a}_i}\in G(M_F)$ then $x^{{\bf a}_i}\in F$ and there exists $\mu\in\K\setminus\{0\}$ such that $x^{{\bf a}_i}-\mu x^{\bf u}\in I$. We have
\[
x^{{\bf a}_i}-\lambda x^{{\bf a}_i+{\bf w}}=(x^{{\bf a}_i}-\mu x^{\bf u})+\mu(x^{\bf u}-\lambda x^{{\bf u}+{\bf w}})-\lambda x^{\bf w}(x^{{\bf a}_i}-\mu x^{\bf u}),
\]   
hence $x^{{\bf a}_i}-\lambda x^{{\bf a}_i+{\bf w}}\in I$ and consequently $x^{{\bf a}_i+{\bf w}}\in F$. Therefore we obtain the desired equality.
\end{proof}

The next example comments on certain subtleties of the above proposition. 

\begin{ex1}
\label{fiber_gen}
{\em Let $I_3$ be the ideal of  Example~\ref{draw_fibers} c). The fiber $F_1=F_{x^3y^3}$ is finite even though $L_{F_1}=\langle (2,-1), (-1, 2) \rangle $ contains  $(1,1)$. Note that $(1,1)\notin L_{pure,F_1}$ and $L_{pure,F_1}=({\bf 0})$. For the fiber $F_2=F_{x^4y^3}$, one can see that  
\[
M_{F_2}=(x^{10}, x^8y,x^6y^2,x^4y^3,x^3y^5,x^2y^7,xy^9,y^{11}),
\] 
thus $L_{pure,F_2}=\langle (1,1), (0,3)\rangle$ and $L^+_{pure,F_2}=(0,3)\N+(1,1)\N+(3,0)\N$. Therefore we have $\{x^{10}x^{b+3c}y^{3a+b}:\ a,b,c\in\N\} \subset F_2$.
}
\end{ex1}

\noindent Let $I$ be a binomial ideal. We define an equivalence relation $"\equiv"$ on the set of $I$-fibers and a partial order  $"<"$ of the equivalence classes which generalize those from \cite{CTV}:

 \begin{def1}\label{equiv}
 {\em If  $F$, $G$ are $I$-fibers, we let $F\equiv G$ if there exist ${\bf u},{\bf v}\in \N^n$ such that $x^{\bf u}F\subset G$ and $x^{\bf v} G\subset F$ and denote by $\overline F$, the equivalence class of $F$. We set $\overline{F}\leq \overline{G}$ if there exists ${\bf u} \in\N^n$ such that $x^{\bf u} F \subset G$. We write $\overline F< \overline G$ if $\overline  F\leq  \overline G$ and $\overline{F}\neq \overline{G}$.}
\end{def1}

\noindent In \cite[Proposition 3.5]{CTV} it was shown for lattice ideals that any two equivalence classes of fibers have the same cardinality. This is no longer necessarily true for an arbitrary binomial ideal $I$ as you can see in the Example~\ref{example_multi} b), but the cardinality of $\overline{F}$ for an $I$-fiber $F$ can be computed similarly by replacing $L_F$ with $L_{pure,F}$. 

\begin{ex1}
\label{example_multi} $ $ \newline
\noindent
{\em a) Let $I_1$ be the ideal of Example~\ref{draw_fibers} a).  The infinite $I_1$-fibers are equivalent.

\noindent b) Let $I_2$ be the ideal of Example~\ref{draw_fibers} b). We note that
\begin{itemize}
\item{} $\overline{F_{y}}=\{F_{y},F_{xy}\}$, 
\item{} $\overline{F_{y^2}}=\{F_{y^2},F_{xy^2}\}$, 
\item{} $\overline{F_{y^3}}=\{F_{y^3}\}$,
\item{} $\overline{F_{y^4}}=\{F_{y^4}\}$. 
\end{itemize}
The set of equivalence classes of $I_2$-fibers is totally ordered, the maximal element being $\overline{F_{y^4}}$ and the minimal one $\overline{F_{1}}$, where $F_{1}=\{1\}$.  

\noindent c) Let $I_3$ be the ideal of Example~\ref{draw_fibers} c). The three infinite fibers $F_{x^4y^3},F_{x^5y^3},F_{x^4y^4} $ are equivalent. The equivalence classes  $\overline{F_{x^6y}}$, $\overline{F_{xy^6}}$ and $\overline{F_{x^3y^3}}$  are incomparable and minimal when restricting to equivalence classes of fibers of cardinality greater than one.}
\end{ex1}

Note that it is possible for a binomial ideal to contain monomials. If this is the case, then it obviously contains infinitely many, but the surprising fact is that all of these monomials form an $I$-fiber, which becomes the maximal $I$-fiber.

\begin{thm1}\label{unique_monomial_fiber}
Let $I$ be a binomial ideal and denote by $F(I)$ the set $\{x^{\bf w}| \ x^{\bf w}\in I\}$. Then $F(I)\neq\emptyset$ if and only if there exist $\lambda\neq\mu\in \K$ and monomials $x^{\bf u},x^{\bf v}$ such that $x^{\bf u}-\lambda x^{\bf v}$ and $x^{\bf u}-\mu x^{\bf v}$ belong to $I$. Furthermore, if $F(I)\neq\emptyset$  then $F(I)$ is an $I$-fiber, $\overline{F(I)}=\{F(I)\}$ and $\overline{F}\leq_I \overline{F(I)}$ for every $I$-fiber $F$.  
\end{thm1}
\begin{proof}
Assume that $F(I)\neq\emptyset$ and let $x^{\bf w}\in I$. If $x^{\bf v}$ is a monomial such that ${\bf v}\sim_I {\bf w}$ then there exists $\lambda\in\K\setminus\{0\}$ such that $x^{\bf u}-\lambda x^{\bf v}\in I$. Since $x^{\bf w}\in I$ then also $x^{\bf v}\in I$, and thus $x^{\bf v}\in F(I)$. Hence $F_{\bf w}\subseteq F(I)$. For the converse inclusion, note that if $x^{\bf v}\in F(I)$ then $x^{\bf v}\in I$ and implicitly $x^{\bf w}-x^{\bf v}\in I$. Therefore ${\bf v}\sim_I {\bf w}$ and $x^{\bf v}\in F_{\bf w}$, which implies $F(I)\subseteq F_{\bf w}$. 

It is easy to see from the definition of $F(I)$ that $x^{\bf w}F(I)\subseteq F(I)$ and $x^{\bf w}F\subseteq F(I)$ for any $I$-fiber $F$. This implies in particular $\overline{F}\leq_I\overline{F(I)}$ for every $I$-fiber $F$. On the other hand, if $x^{\bf v}F(I)\subseteq F$ for some $I$-fiber $F$ with $F\neq F(I)$ then we obtain a contradiction from Theorem~\ref{partition_classif}, since $x^{\bf v}F(I)\subseteq F(I)$. Therefore it follows that $\overline{F(I)}=\{F(I)\}$, and we are done.
\end{proof}

\begin{ex1}\label{examples_monomial_fiber}
\rm Let $I_1,\ldots,I_4$ be the ideals of Example~\ref{draw_fibers}. Then: 
\begin{enumerate}
\item[(a)] $F(I_1)=\emptyset$;
\item[(b)] $F(I_2)=\emptyset$ if $\chara(\K)\neq 3$, and $F(I_2)=F_{y^4}$ when $\chara(\K)=3$; 
\item[(c)] $F(I_3)=\emptyset$;
\item[(d)] $F(I_4)=F_{x^4y^3}\cup F_{x^3y^4}\cup F_{x^4y^4}$ if $\chara(\K)\neq 5$, where $F_{x^4y^3}, F_{x^3y^4}, F_{x^4y^4}$ are the underlying sets of the  corresponding fibers of $I_3$; and $F(I_4)=F_{x^6y}$ if $\chara(\K)=5$, where $F_{x^6y}$ is the corresponding fiber of $I_4$.  
\end{enumerate}
\end{ex1}

\begin{rem1}\label{construction_monomial_fiber}
\rm Whether a binomial ideal $I$ contains or not monomials is computationally detectable through a single Gr\" obner basis computation. Indeed, if we compute a reduced Gr\" obner basis $\mathcal{G}$ with respect to any monomial order then $\mathcal{G}$ contains a monomial if and only if $I$ contains monomials. The proof of this remark follows immediately from Buchberger's algorithm for computing a reduced Gr\" obner basis. 
\end{rem1}

The proof of \cite[Theorem 3.8]{CTV} applies ad litteram for an arbitrary binomial ideal $I$ and is based on the noetherian property of a chain of monomial ideals associated to the fibers.  

\begin{thm1} Let $I$ be a binomial ideal. Any descending chain of equivalence classes of fibers
\[
\overline F_{1} \ > \ \cdots\  >{ }\   \overline F_k\ > \ \overline F_{k+1}\  > \  \cdots
\]
is finite.  
\end{thm1}

It is easy to see that if $F$ is finite then $\overline{F}=\{F\}$. In the next example we compute the equivalence classes of the infinite fibers for the binomial ideals of Example~\ref{draw_fibers}.

Let $F$ be a fiber of a binomial ideal $I$. We define two binomial ideals contained in $I$
\[
I_{<\overline F}=(B\in I:\ B\ \text{ binomial},\ \overline F_B < \overline F),  \ \ \  I_{\leq\overline F}=(B\in I:\ \overline F_B\leq  \overline F).
\]
We say that $F$ is a {\it Markov fiber} if there exists a minimal binomial generating set $S$ for $I$ and $B\in S$ such that $F_B\equiv F$. The definition of the equivalence relation among the fibers together with the induced partial order allows us to identify the equivalent classes of the Markov fibers and to prove that the set of equivalence classes of Markov fibers is an invariant of the ideal. This is the content of the next theorem. 

\begin{thm1}\label{markov_fibers} Let $I$ be a binomial ideal, $S$ a binomial generating set of $I$ and $F$ an $I$-fiber. Then   
\[
I_{<\overline F}=(B:\  B\in S,\  \overline{F_{B}}<  \overline F) \text{ and } I_{\leq \overline F}=(B:\  B\in S,\  \overline{F_{B}}\leq  \overline{F})\ . 
\]  Thus $F$ is a Markov fiber of $I$ if and only if $I_{< \overline F}\neq I_{\leq \overline F}$. Moreover the set $\{\overline{F_{B}}:\  B\in S\}$ is an invariant of $I$.
\end{thm1}
\begin{proof} We will show the statement for $I_{<\overline F}$, the other one having a similar proof. Let $J=(B:\  B\in S,\  \overline F_{B}<_I \overline F)$. It is clear that $J\subset I_{< \overline F}$. For the other inclusion it is enough to prove that if $B=x^u-\lambda x^v\in I_{<\overline F}$ then $B\in J$. Let $B=x^u-\lambda x^v\in I_{<\overline F}$. Since $B\in I$ then $B=\sum_{j=1}^t c_j x^{a_{j}}B_{j}$, where $B_{j}=x^{u_j}-\lambda_j x^{v_j}\in S$ are not necessarily distinct. The expression $ B=\sum_{j=1}^t c_j x^{a_{j}}B_{j} $ implies that, after a possible rearrangement of indices, there exist binomials $B_1,\ldots, B_s$ and a sequence of monomials $x^{m_1},\dots , x^{m_{s+1}}$ in the fiber $F_B$ such that $x^{m_1}=x^u=x^{b_1}x^{u_1}, \dots, x^{m_i}=x^{b_{i}}x^{v_{i}}=x^{b_{i+1}}x^{u_{i+1}}, \dots, x^{m_s}=x^{b_{s}}x^{v_{s}}=x^v$. Then $x^u-\lambda_1 \dots \lambda_{s}x^v=\sum_{i=1}^s \lambda_1 \dots \lambda_{i-1} x^{b_{i}}B_{i} $. Note also that since all monomials $x^{m_1},\dots , x^{m_{s+1}}$ belong to  the fiber $F_B$ it follows that $x^{b_{i}}  F_{B_i}\subset F_B$. Thus $\overline F_{B_i}\le_I \overline F_B$. Since $\overline  F_B<_I\overline  F$ we see that $\overline F_{B_i}<_I \overline F$.
There are two cases to analyze. In the first case $\lambda=\lambda_1 \dots \lambda_{s}$, and we have $B\in J$. In the second case, $\lambda\neq\lambda_1 \dots \lambda_{s}$ implies $x^u-\lambda x^v\in I$ and $x^u-\lambda_1\dots \lambda_{s}x^v\in I$, and by Theorem~\ref{unique_monomial_fiber} we have  $F(I)\neq\emptyset$ and $F_B=F(I)$. Again by Theorem~\ref{unique_monomial_fiber} we obtain a contradiction to $\overline F_{B}<_I \overline F$, and we are done.
\end{proof}

For a lattice ideal $I$ we have even more: if $S_1$, $S_2$ are two minimal binomial generating sets of $I$ of minimal cardinality, then $\{\overline{F_{B}}:\  B\in S_1\}= \{\overline{F_{B}}:\  B\in S_2\}$ where the equality holds for the {\it multisets} (i.e. sets together with the multiplicities of their elements), see \cite[Corollary 4.14]{CTV}. 

\begin{ex1}[]
{\em  Let $I_2$ be the ideal of  Example~\ref{draw_fibers} b), which is minimally generated by the four binomials. The set of equivalence classes of Markov fibers of $I_2$ is the set $\{ \overline{F_y}, \overline{F_{y^3}}, \overline{F_{y^4}}\}$. It is not hard to see that $I_2$ has also a minimal generating set of cardinality three: $I_2=(y -x^2y,y^3-xy^3,y^4-y^5)$.
}
\end{ex1}

\section{Indispensable Monomials and Binomials}

Let $I$ be a binomial ideal and $F$ an $I$-fiber.  We consider the monomial ideals $M_{\overline F}$ generated by the monomials of {\em all} fibers $G$ with $G\equiv F$, and $M_I$ generated by all monomials $x^{\bf u}\in \supp(B)$, where $B\in I$ is a nonzero binomial. It is clear that $M_F\subseteq M_{\overline F}$.  We  note the following:

\begin{lem1}\label{generating_M_I}
If $I=(x^{{\bf u}_1}- \lambda_1 x^{{\bf v}_1},\ldots,x^{{\bf u}_s}- \lambda_s x^{{\bf v}_s})$ with $\lambda_1,\ldots,\lambda_s\in\K\setminus\{0\}$ then $M_I=(x^{{\bf u}_1},x^{{\bf v}_1},\ldots,x^{{\bf u}_s},x^{{\bf v}_s})$.
\end{lem1}
\begin{proof}
One inclusion is immediate. For  "$\subseteq$" suppose that $x^{\bf u}- \lambda x^{\bf v}\in I$ for some $\lambda\in\K\setminus\{0\}$. Hence there are  polynomials $f_j\in R$ such that  
\[
x^{\bf u}- \lambda x^{\bf v}=\sum_{j=1}^s f_j(x^{{\bf u}_j}- \lambda_jx^{{\bf v}_j})\ .
\]
Thus for some $i, k\in [s]$, $x^{{\bf u}_i}$ or $x^{{\bf v}_i}$ divides $x^{\bf u}$, and  $x^{{\bf u}_k}$ or $x^{{\bf v}_k}$ divides $x^{\bf v}$. Therefore $x^{\bf u}, x^{\bf v}$ belong to $(x^{{\bf u}_1},x^{{\bf v}_1},\ldots,x^{{\bf u}_s},x^{{\bf v}_s})$.
\end{proof}

We  also note that if $x^{\bf u}$ is a minimal monomial generator of $M_I$ then $x^{\bf u}$ is also a minimal monomial generator of $M_{\overline{F_{\bf u}}}$ and of $M_{F_{\bf u}}$.

\begin{ex1}[]
\label{monomial_fiber}
{\em Let $I_2$ be the ideal of Example~\ref{example_multi} b). Then $M_{I_2}=(y)$. Also, if $F=F_{xy^2}$ then $M_{F}=(xy^2)$ and $M_{\overline {F}}=(y^2)$. If   $G=F_{y}$ then $M_{G}=M_{\overline{G}}=M_{I_2}$.}
\end{ex1}

In \cite[Proposition 1.5]{KO} it was shown that if $I$ is an $A$-homogeneous binomial ideal for some $A\subset \Z^d$ such that  $\N A\cap (-\N A)=\{{\bf 0}\}$, then the indispensable monomials of $I$ are exactly the minimal monomial generators of $M_I$. The same proof applies to any binomial ideal $I$.  We isolate this remark:

\begin{rem1}\label{known_indisp}
{\em Let $I$ be a binomial ideal and $S$ a system of binomial generators  of $I$. The indispensable monomials of $I$ are precisely the elements of  $G(M_I)$. Moreover  $G(M_I)$ comes from $\bigcup_{B\in S} \supp(B)$ by   keeping the minimal elements according to divisibility.}
\end{rem1}

Next, for $F$ an $I$-fiber, we define the  graph $\Gamma_F$ (on possibly infinitely many vertices): 

\begin{def1}\label{gamma_graph}
{\em  $\Gamma_F$ is the graph with vertices the elements of  $F$ and edges $\{x^{\bf u},x^{\bf v}\}$ whenever $x^{\bf u}-\lambda x^{\bf v}\in I_{<\overline{F}}$ for some $\lambda\in\K\setminus\{0\}$.
}
\end{def1}

The connected components of $\Gamma_F$ determine the indispensable monomials and binomials as the next two theorems show.

\begin{thm1}\label{indisp_mon_fiber}
Let $I$ be a binomial ideal. The monomial
 $x^{\bf u}$ is a minimal monomial generator of $M_I$ (and thus an indispensable monomial) if and only if the following conditions hold simultaneously
\begin{enumerate}
\item[(a)] $|F_{\bf u}|\geq 2$, 
\item[(b)] $x^{\bf u}$ is an isolated vertex of $\Gamma_{F_{{\bf u}}}$,
\item[(c)] $x^{\bf u}$ is a minimal generator of $M_{\overline{F_{\bf u}}}$.
\end{enumerate}
\end{thm1}

\begin{proof}
Assume first that $x^{\bf u}$ is a minimal monomial generator of $M_I$. Conditions (a) and (c) follow immediately. For condition (b) we suppose that $x^{\bf u}$ is not an isolated vertex of $\Gamma_{F_{\bf u}}$: there exists $x^{\bf v}\in F_{\bf u}$ so that $\{x^{\bf u},x^{\bf v}\}$ is an edge of $\Gamma_{F_{\bf u}}$ and thus $x^{\bf u}-\lambda x^{\bf v}\in I_{<\overline{F_{\bf u}}}$ for some $\lambda\in\K\setminus\{0\}$. By Theorem~\ref{markov_fibers}, the binomial $B=x^{\bf u}-\lambda x^{\bf v}$ is not part of any minimal binomial generating set of $I$. On the other hand, since $x^{\bf u}\in G(M_I)$, it follows by Lemma~\ref{generating_M_I} that for any minimal binomial generating set $S$ there exists a binomial $x^{\bf u}-\mu x^{\bf w}\in S$ with $\mu\in\K\setminus\{0\}$. Fix such a minimal system of binomial generators $S$ and let $B_i=x^{\bf u}-\lambda_i x^{{\bf w}_i}$ with $i=1,\ldots,t$ be the binomials in $S$ having in their support $x^{\bf u}$. If we define $B'_i=x^{\bf v}-(\lambda_i/\lambda) x^{{\bf w}_i}$ for $i=1,\ldots,t$, then one can easily see that the set $S'=(S\setminus \{B_1,\ldots,B_t\})\cup\{B,B'_1,\ldots,B'_t\}$ is a system of generators of $I$. Minimizing $S'$ we obtain a minimal system of binomial generators $S''$ which is properly contained in $S'\setminus\{B\}$ since $B$ is not part of any minimal system of binomial generators. Thus $x^{\bf u}$ does not appear in the support of any binomial from $S''$, and therefore by Lemma~\ref{generating_M_I}, $x^{\bf u}$ is not a minimal monomial generator of $M_I$, a contradiction.

For the converse, assume that  $x^{\bf u}$ satisfies conditions (a)-(c). Suppose that $ x^{\bf u}$ is not a minimal monomial generator of $M_I$. Since $|F_{\bf u}|\geq 2$, we conclude that there exists a monomial $x^{{\bf v}}$ such that $x^{{\bf u}}-{\lambda}x^{{\bf v}}\in I$ for some $\lambda\in\K\setminus\{0\}$. Since $ x^{\bf u}\notin G(M_I)$ we conclude that there exists a minimal binomial system of generators $S$ of $I$ and $B\in S$ so that $x^{\bf u}=x^{\bf w}x^{\bf u'}  $with $x^{\bf u'}\in\supp(B)$  and $x^{\bf w}\neq 1$. In particular we get that  $\overline{F_{{\bf u'}}}\leq\overline{F_{\bf u}}$. We first examine the case when  $B\in I_{<\overline{F_{\bf u}}}$. If $\supp(B)=\{ x^{{\bf u'}}, x^{{\bf v'}}\}$ we let $x^{{\bf w'}}=x^{{\bf v'}}\cdot x^{\bf w}\in F_{\bf u}$. This implies that $\{x^{{\bf u}},x^{{\bf w'}}\}$ is an edge of $\Gamma_{F_{\bf u}}$,  contradicting condition (b). If $B\notin I_{<\overline{F_{\bf u}}}$ then {  $F_{{\bf u}}\equiv F_{{\bf u}'}$} and since $x^{{\bf u'}}$ properly divides $x^{{\bf u}}$ it follows that  $x^{{\bf u}}$ is not a minimal generator of $M_{\overline{F_{\bf u}}}$, a contradiction to (c).
\end{proof}

The description of the indispensable monomials in terms of the graph of their fibers will be useful in the next theorem.

\begin{thm1}\label{indisp_bin_fiber}
Let $I$ be a binomial ideal. The binomial $B\in I$ is indispensable if and only if the graph $\Gamma_{F_{B}}$ consists of two isolated vertices.
\end{thm1} 
\begin{proof}
Let $B=x^{\bf u}-\lambda x^{\bf v}$ for some $\lambda\in\K\setminus\{0\}$ and $F=F_B$. Assume first that the graph $\Gamma_{F}$ consists of two isolated vertices, that is $F=\supp(B)$ and $B\notin I_{<\overline{F}}$. By Theorem~\ref{markov_fibers}, any generating set of $I$ contains a binomial $B'$ such that $F_{B'}\equiv F$. However since $F$ is finite it follows that $\overline{F}=\{F\}$ and since $F=\supp(B)$ it must be that up to a multiplication by a nonzero coefficient $B'=x^{\bf u}-\mu x^{\bf v}$ for some $\mu\in\K\setminus\{0\}$. If $\mu\neq\lambda$ then $x^{\bf u},x^{\bf v}\in I$, thus $F=F(I)$ is an infinite fiber, a contradiction to the finiteness of $F$. Hence $\mu=\lambda$ and so $B$ is indispensable. 

For the other direction, assume that $B$ is indispensable. Thus the elements of $\supp(B)$ are indispensable monomials. By Theorem~\ref{indisp_mon_fiber} the elements of  $\supp(B)$ are isolated vertices of $\Gamma_{F}$ and $B\notin I_{<\overline{F}}$. It remains to show that $F=\supp(B)$. Suppose that
$x^{\bf w}\notin \supp(B)$ is also in $F$. This implies that there exist $\lambda_1,\lambda_2\in\K\setminus\{0\}$ such that $x^{\bf u}-\lambda_1x^{\bf w},x^{\bf v}-\lambda_2 x^{\bf w}\in I$. Hence also $x^{\bf u}-(\lambda_1/\lambda_2)x^{\bf v}$ belongs to $I$. Let $S$ be any minimal generating set of $I$. It necessarily contains $B$. We have two cases to analyze. First assume that $\lambda_1/\lambda_2=\lambda$. This implies that set $S\setminus\{B\}\cup\{x^{\bf u}-\lambda_1x^{\bf w},x^{\bf v}-\lambda_2x^{\bf w}\}$  is a system of generators of $I$ and when we minimize it we obtain a minimal system of generators of $I$ not containing $B$, a contradiction since $B$ is indispensable. In the second case $\lambda_1/\lambda_2\neq\lambda$, and we obtain $x^{\bf u},x^{\bf v}\in I$. Therefore $F=F(I)$, $x^{\bf w}\in I$ and it can be seen immediately that $S\setminus\{B\}\cup\{x^{\bf u}-\lambda x^{\bf w},x^{\bf v}-x^{\bf w}\}$ is a system of generators of $I$. Minimizing it we obtain a minimal system of generators not containing $B$, a contradiction since $B$ is indispensable. Hence the graph $\Gamma_{F}$ consists of two isolated vertices $x^{\bf u},x^{\bf v}$.    
\end{proof}

\begin{ex1}
\label{indisp_mon_binom}
{\em Let $I_3$ be the ideal of Example~\ref{draw_fibers} c). Since  
\[
G(M_{I_3})=\{x^8,x^6y,x^5y^2,x^3y^3,x^2y^5,xy^6,y^8\}
\] 
then $I_3$ has seven indispensable monomials by Remark~\ref{known_indisp}. This implies that a minimal generating set of $I_3$ can have cardinality no less than four. The indispensable binomials of $I_3$ are $x^6y-x^8$ and $xy^6-y^8$ as follows from Theorem~\ref{indisp_bin_fiber} and the study of the fibers of the indispensable monomials. See Examples~\ref{draw_fibers} c), \ref{example_multi} c). }
\end{ex1}

The following result generalizes \cite[Corollary 1.11]{KO}. 

\begin{cor1}\label{indisp_subideal}
Let $J\subset I$ be two binomial ideals. If $B$ is an indispensable binomial of $I$ and $B\in J$ then $B$ is indispensable in $J$. 
\end{cor1}
\begin{proof}
By Theorem~\ref{indisp_bin_fiber} we have that the $I$-fiber $F_B$ is equal to $\supp(B)$. Since the fiber of $B$ in $J$ is a subset of $F_B$ and contains  $\supp(B)$, then it equals $F_B$. Moreover since $J_{<\overline{F_B}}\subset I_{<\overline{F_B}}$ it follows that $B\notin J_{<\overline{F_B}}$. The conclusion now is obtained by applying Theorem~\ref{indisp_bin_fiber} one more time.
\end{proof}

\section{Computing  indispensable  binomials of  $I$}

Let $I$ be a binomial ideal and $S$ a system of binomial generators of $I$. 
In this section we provide an algorithm for finding the indispensable binomials of $I$.  This algorithm   generalizes by far the three algorithms known to us which are given in the restrictive case of positively $A$-graded toric ideals (i.e. $\N A\cap (-\N A)=\{{\bf 0}\}$), which are in particular pure binomial ideals.
We recall that:
\begin{itemize}
\item the algorithm in  \cite[Theorem 2.4]{OH1}  implies computation of $n!$ reduced Gr{\" o}bner bases with respect to  the lexicographic orders, 
\item the algorithm in \cite[Theorem 3.4]{CKT}  implies the computation of one Gr{\" o}bner basis   and the knowledge of the minimal elements in the set of $I$-fibers,
\item the algorithm in  \cite[Theorem 13]{OV} implies computation of  $n$ reduced Gr{\" o}bner basis with respect to $n$ degree reverse lexicographic orders.
\end{itemize}
The simplicity and the swiftness of the algorithm we propose in this section depends only on the information given by $S$. First we define  a graph with vertices the (at most 2$|S|$) elements in the union of the supports of the binomials of $S$.

\begin{def1}
{\em Let $\mathcal F(S)$ be the graph whose vertices are the monomials in $\bigcup_{B\in S} \supp(B)$ and  edges $\{x^{\bf u},x^{\bf v}\}$ whenever, up to a nonzero scalar multiplication, $x^{\bf u}-\lambda x^{\bf v}\in S$ for some $\lambda\in\K\setminus\{0\}$.} 
\end{def1}

Note that the graph $\mathcal F(S)$ may not be simple, since it may have multiple edges as the following example shows. 

\begin{ex1}
\label{forest}
{\em Let $I\subset\C[x,y]$ be the binomial ideal generated by $S=\{y^8-xy^6,x^2y^5-x^3y^3,x^3y^3-x^5y^2,x^6y-x^8,x^6y-2x^8\}$. The graph $\mathcal F(S)$ is   depicted below:

\begin{figure}[hbt]
\begin{center}
\psset{unit=1cm}
\begin{pspicture}(3.75,1.7)(5,3.8)
\rput(0,3){$\bullet$}
\rput(1.5,3){$\bullet$}
\rput(0,2.6){$y^8$}
\rput(1.5,2.6){$xy^6$}
\psline[linewidth=0.7pt](0,3)(1.5,3)
\rput(3.5,3){$\bullet$}
\rput(4.5,2.5){$\bullet$}
\rput(4.5,3.5){$\bullet$}
\rput(3.3,2.6){$x^3y^3$}
\rput(4.7,2.25){$x^2y^5$}
\rput(4.7,3.75){$x^5y^2$}
\rput(4,1.5){$\mathcal F(S)$}
\psline[linewidth=0.7pt](4.5,2.5)(3.5,3)
\psline[linewidth=0.7pt](3.5,3)(4.5,3.5)
\rput(6.5,3){$\bullet$}
\rput(8,3){$\bullet$}
\rput(6.5,2.6){$x^8$}
\rput(8,2.6){$x^6y$}
\pscurve[linewidth=0.7pt](6.5,3)(7.25,3.2)(8,3)
\pscurve[linewidth=0.7pt](6.5,3)(7.25,2.8)(8,3)
\end{pspicture}
\end{center}
\caption{The graph associated to a generating set $S$ of a binomial ideal}
\label{Fig2}
\end{figure}
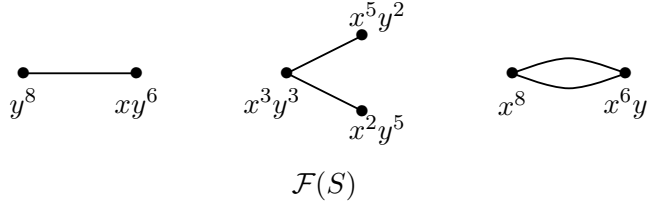

}

\end{ex1}

Moreover, if $S$ is a minimal binomial generating set for a binomial ideal $I$ with $F(I)=\emptyset$ then $\mathcal F(S)$ is a forest.

\begin{thm1}\label{algorithm_indisp_thm}
Let $S$ be a binomial generating set of $I$. The binomial  $B\in S$ is   indispensable  if and only if $\supp(B)\subset G(M_I)$ and the induced graph on the vertices of $\supp(B)$ is a connected component of $\mathcal F(S)$ consisting of a simple edge.
\end{thm1}
\begin{proof}
If $B\in S$ is indispensable then the first assertion follows from Remark~\ref{known_indisp}. Moreover by Theorem~\ref{indisp_bin_fiber}, $|F_B|=2$. Thus  
 the induced graph on the vertices of $\supp(B)$ is necessarily a connected component of $\mathcal F(S)$. Furthermore, if this induced graph was not a simple edge then there would exist $\lambda_1\neq\lambda_2\in\K\setminus\{0\}$ such that $x^{\bf u}-\lambda_1 x^{\bf v},x^{\bf u}-\lambda_2 x^{\bf v}\in S$. This implies that $x^{\bf u},x^{\bf v}\in I$, thus $F_B=F(I)$ is an infinite fiber, a contradiction since $|F_B|=2$.

For the converse, assume that $B=x^{\bf u}-\lambda x^{\bf v}$ for some $\lambda\in\K\setminus\{0\}$, $\supp(B)\subset G(M_I)$ and the induced graph on the vertices of $\supp(B)$ is a connected component of $\mathcal F(S)$ consisting of a simple edge. The last condition implies that $\supp(B)\cap \supp(B')=\emptyset$  for all $B'\in S$, $B'\neq B$. We can assume that $S=\{ B_1, B_2,\ldots, B_s\}$ (where $B=B_1$) and that $B_j=x^{{\bf u}_j}-\lambda_j x^{{\bf v}_j}$. Suppose now that there exists a monomial $x^{{\bf w}}\in F_{B}\setminus \supp(B)$. This implies that there exists $\mu\in\K\setminus\{0\}$ such that $x^{{\bf u} }- \mu x^{{\bf w}}\in I$. Thus
\begin{eqnarray}\label{one}
x^{{\bf u} }- \mu x^{{\bf w}}=\sum_{j,t} c_{\alpha_{j,t}}x^{\alpha_{j,t}}B_j\ ,
\end{eqnarray}
where $c_{\alpha_{j,t}}\in\K$ and the monomials $x^{\alpha_{j,t}}$ are such that $x^{\alpha_{j,t_1}}\neq x^{\alpha_{j,t_2}}$ for any $j$ and $t_1\neq t_2$. 
Since $x^{\bf u}\in G(M_I)$  and $x^{\bf u}\notin \supp(B_j)$ for $j\neq 1$ we can assume  that  $c_{\alpha_{1,1}}=1$ and $x^{\alpha_{1,1}}=1$.
By our assumption $x^{{\bf w}}\neq x^{{\bf v}}=x^{{\bf v}_1}$. Thus  $ x^{{\bf v}_1}$ must appear at least twice in the RHS   of Equation  \ref{one}. It follows that  $ x^{{\bf v}_1}$ is divisible by  $x^{{\bf u}_j}$ or $x^{{\bf v}_j}$ for some $j\neq 1$. Since $x^{{\bf v}_1}\in G(M_I)$   
this would imply that  $x^{{\bf v}_1}=x^{{\bf u}_j}$ or $x^{{\bf v}_1}=x^{{\bf v}_j}$. However this is impossible since the induced graph on the vertices of $\supp(B)$ is a connected component of $\mathcal F(S)$. Therefore our assumption is wrong and $F_B=\supp(B)$. 
Assume now that  $B\in I_{<\overline{F_B}}$.  By Theorem~\ref{markov_fibers} 
\begin{eqnarray}\label{two}
x^{{\bf u} }-\lambda x^{{\bf v}}=\sum_{i,s} c_{\beta_{i,s}}x^{\beta_{i,s}}B_i\ ,
\end{eqnarray}
where the binomials $B_i$ that appear in the RHS have the property that  $\overline{F_{B_i}}<\overline{F_B}$. Thus   $x^{{\bf u}}$ is properly divisible by some monomial in $\supp(B_i)$ for  $i\neq 1$. Since $x^{{\bf u}}\in G(M_I)$ this is a contradiction. Thus $B\notin I_{<\overline{F_B}}$ and $\Gamma_{F_B}$ consists of two isolated vertices. The desired conclusion follows from Theorem~\ref{indisp_bin_fiber}.
\end{proof}

\begin{ex1}\label{indispensable_examples}
\rm Note that for the ideal $I$ of Example~\ref{forest} there is only one connected component consisting of a simple edge, namely the one corresponding to $\supp(y^8-xy^6)$. Thus $I$ has just one indispensable binomial.
\end{ex1}

The following algorithm is an immediate application of Theorem~\ref{algorithm_indisp_thm}.

\begin{algorithm}[h] 
\caption{Computing the indispensable binomials  of a binomial ideal $I$}  
\label{indisp_binomials}

\begin{algorithmic}[1]
\REQUIRE $F = \{B_1,\ldots, B_s\}\subseteq \K[X]$, with $B_i=x^{{\bf u}_i}- \lambda_i x^{{\bf v}_i}$ for $i\in [s]$, where $\lambda_1,\ldots,\lambda_s\in\K\setminus\{0\}$.
\ENSURE $F'\subset F$, the set of indispensable binomials of $I=(B_1,\ldots, B_s)$.
\vspace{0.1cm}
\STATE Compute $G(M_I)$, a subset of $\{x^{{\bf u}_1},x^{{\bf v}_1},\ldots,x^{{\bf u}_s},x^{{\bf v}_s}\}$, and set $T=\{i: \{x^{{\bf u}_i},x^{{\bf v}_i}\}\subset  G(M_I)\}$.
\STATE If {$T=\emptyset$} then $F'=\emptyset$.
\STATE Otherwise, for every $i\in T$ check whether $x^{{\bf u}_i}\in\Supp(B_j)$ or $x^{{\bf v}_i}\in\Supp(B_j)$ for some $j\neq i$. 
\STATE $F'=\{x^{{\bf u}_i}-\lambda_i x^{{\bf v}_i}: i\in T, x^{{\bf u}_i}\notin\Supp(B_j),x^{{\bf v}_i}\notin\Supp(B_j), \text{ for all } j\neq i\}$.
\end{algorithmic}
\end{algorithm}

\begin{rem1}\label{polynomial_time_complexity}
\rm We note that step 1 of Algorithm~\ref{indisp_binomials} involves checking the divisibility relations of the elements of  $\{x^{{\bf u}_1},x^{{\bf v}_1},\ldots,x^{{\bf u}_s},x^{{\bf v}_s}\}$. Checking one such divisibility is equivalent to computing the difference vector of the two exponent vectors and see if it belongs to $\N^n$. Thus the running time of checking such a divisibility is $O(n)$. Since there are $\frac{2s(2s-1)}{2}$ such divisibility relations to check, the total running time of having the output of step 1 is $O(s^2n)$. Therefore Algorithm~\ref{indisp_binomials} is a polynomial-time algorithm. 
\end{rem1}

\begin{cor1}\label{gen_indispensable_cor}
Let $I$ be a binomial ideal minimally generated by $s$ binomials. $I$ is generated by indispensable binomials if and only if $|G(M_I)|=2s$. 
\end{cor1}
\begin{proof}
Assume that $I$ is minimally generated by the binomials $x^{{\bf u}_1}-\lambda_1 x^{{\bf v}_1}$, $\ldots$, $x^{{\bf u}_s}-\lambda_s x^{{\bf v}_s}$, where $\lambda_1,\ldots,\lambda_s\in\K\setminus\{0\}$. If the binomials are indispensable then $|G(M_I)|=2s$ by  Theorem~\ref{algorithm_indisp_thm}. Conversely, if $|G(M_I)|=2s$ we obtain the desired conclusion from Algorithm~\ref{indisp_binomials}. 
\end{proof}

\begin{algorithm}[h] 
\caption{Testing whether a binomial ideal is generated by indispensable binomials}  
\label{gen_indispensable}

\begin{algorithmic}[1]
\REQUIRE $F = \{B_1,\ldots, B_s\}\subseteq \K[X]$, a set of binomials generating $I$.
\ENSURE Is $I$ generated by indispensable binomials? YES or NO
\vspace{0.1cm}
\STATE Compute $S\subset F$, a set of minimal generators for $I$.
\STATE Compute $G(M_I)$ from $S$.
\STATE If $|G(M_I)|=2|S|$ then $I$ is generated by indispensable binomials, otherwise not. 
\end{algorithmic}
\end{algorithm}

As an example we  recover immediately that the toric ideal of $\mathcal{A}_{333}$, the $3\times 3\times 3$ contingency tables having fixed two-dimensional marginal totals, is generated by indispensable binomials, \cite[Theorem 1]{AT}. First, using  CoCoA \cite{Co} we compute  a generating set for $I_{\mathcal A_{333}}$: it has cardinality $114$. Then  we minimize this set to get a generating set of cardinality $81$ and finally we compute $|G(M_{I_{\mathcal A_{333}}})|$, which is 162. The criterion of  Corollary~\ref{gen_indispensable_cor} finishes the proof.

\section{Graver bases of binomial ideals}

\medskip

The notion of primitive binomials plays an important role in the theory of lattice ideals and all applications in terms of computations. We generalize this notion to arbitrary pure binomial ideals, see also \cite{DuZi}.

\begin{def1}\label{graver_bin}
{\em Let $I$ be a  pure binomial ideal. A binomial $0\neq x^{\bf u}-x^{\bf v}\in I$ is called a {\em primitive binomial} of $I$ if there exists no other binomial  $0\neq x^{\bf u'}-x^{\bf v'}\in I$ such that $x^{\bf u'}$ divides $x^{\bf u}$ and $x^{\bf v'}$ divides $x^{\bf v}$. 
The set of all primitive binomials of $I$ is called the {\em Graver basis} of $I$, and denoted by $\Gr(I)$.}
\end{def1}


As in the case of lattice ideals, one can generalize \cite[Lemma 4.6]{St} to show that all elements of the universal 
Gr\"obner basis of $I$ are primitive.  Below
we include the proof for completeness. 

\begin{prop1}\label{ugb_gr} 
Let $I$ be a pure binomial ideal. Every binomial in  the universal Gr\" obner basis of $I$ is contained in $\Gr(I)$. In particular, $\Gr(I)$ is a generating set for the ideal $I$. 
\end{prop1}
\begin{proof}
We argue by contradiction. Assume that $f=x^{\bf u}-x^{\bf v}\notin \Gr(I)$ and $f\in \mathcal{G}_{<}$, the reduced Gr\" obner basis  of $I$ 
according to  the monomial order "$<$". Moreover suppose that  $\ini_{<}(f)= x^{\bf u}>x^{\bf v}$, i.e.   $x^{\bf u}\in G(\ini_<(I))$.
Since $f$ is not primitive,  there exists   $g=x^{\bf u'}-x^{\bf v'}\in I$ such that $f\neq g$ and  $x^{\bf u'}|x^{\bf u}$, $x^{\bf v'}|x^{\bf v}$.   If  $\ini_<(g)=x^{\bf v'}$ then $x^{\bf v'}\in\ini_{<}(I)$ and thus  $x^{\bf v}$ is divisible by an element of 
$G(\ini_<(I))$. This leads to a contradiction since $f$ belongs to the reduced Gr\" obner basis $\mathcal{G}_{<}$ of $I$.  
If  $\ini_{<}(g)=x^{\bf u'}$   then   $x^{\bf u'}=x^{\bf u}$. Thus
 $f-g=x^{\bf v}-x^{\bf v'}\in I$ and $\ini_<(x^{\bf v}-x^{\bf v'})=x^{\bf v}$. This, as before, leads to a contradiction. 
\end{proof}

To prove that $Gr(I)$ is finite, first we remark that if $S$ is an infinite set of monomials then an arbitrary element of $S$ is divisible by some  element of   $G(\langle S\rangle)$.  Thus if  whenever $m, m' \in S$ we have that $m$ does not divide $m'$ neither $m'$ divides $m$ then   $S$ is necessarily finite.   

\begin{prop1}\label{finit_grav} Let $I$ be a pure binomial ideal. The Graver basis of $I$, $\Gr(I)$, is a finite set.
\end{prop1}

\begin{proof} We consider the set  $S=\{ x^{{\bf u}}y^{{\bf v}},\ x^{{\bf v}}y^{{\bf u}}:\   x^{\bf u}-x^{\bf v}\in \Gr(I)\}$.  It is immediate that there are no divisibility relations among distinct elements of $S$. Thus $S$ is finite. It follows that $\Gr(I)$ is finite.
\end{proof}

One can think of the pairs that form $S$ in the above proof, as the support of   binomials   that generate an ideal closely resembling the binomial ideal of the Lawrence lifting of $I$,  see \cite[Theorem 7.1]{St}.  
 In the next theorem we study this ideal. 

\begin{thm1}\label{arbitrary_indispensable}
Let $I\neq (0)$ be a pure binomial ideal and let
 $$\Lambda(I):=(x^{{\bf u}}y^{{\bf v}}-x^{{\bf v}}y^{{\bf u}}:\   x^{{\bf u}}-x^{{\bf v}}\in\Gr(I)) \subset \K[x_1,\ldots,x_n,y_1,\ldots,y_n].$$ 
The set $\{x^{{\bf u}}y^{{\bf v}}-x^{{\bf v}}y^{{\bf u}}: \  x^{{\bf u}}-x^{{\bf v}}\in\Gr(I)\}$ is a minimal system of generators of $\Lambda(I)$ consisting of indispensable binomials.  
\end{thm1} 
\begin{proof}
The conclusion follows immediately by  Algorithm~\ref{gen_indispensable} and Proposition \ref{finit_grav}.
\end{proof}

\begin{rem1}\label{sub_indisp}
{\em Note that the conclusion of Theorem \ref{arbitrary_indispensable} holds for any binomial ideal $J$, where $J=(x^{\bf u}y^{\bf v}-x^{\bf v}y^{\bf u}:\   x^{\bf u}-x^{\bf v}\in A, \ \emptyset\neq A\subset\Gr(I))$.} 
\end{rem1}

We remark that when $I=I_A$, a toric ideal, the ideal $\Lambda(I)$ is equal to the toric ideal of the second Lawrence lifting of $A$, see \cite[Theorem 7.1]{St}. In that setting the Graver basis and the universal Gr\"obner basis of $\Lambda(I)$ coincide, and consist of indispensable binomials. However, in the general case examined in Theorem~\ref{arbitrary_indispensable} the minimal generating set of $\Lambda(I)$ is not necessarily a Graver basis of $\Lambda(I)$ or the universal Gr\"obner basis of $\Lambda(I)$.  
 
\begin{ex1}{\em Let $I_3$ be the ideal of Example \ref{draw_fibers}. Then  $\Gr( I_3)$ has 35 elements, while   the Graver basis of $\Lambda(I_3)$ consists of at least 113 elements. Moreover the reduced Gr\"obner bases of $\Lambda(I_3)$ with respect to lex and degrevlex order differ.  }
\end{ex1}

\end{document}